\documentstyle[12pt]{article}
\input{amssymb.sty}
\newtheorem{theorem}{Theorem}[section]
\newtheorem{prop}[theorem]{Proposition}
\newtheorem{defn}[theorem]{Definition}
\newtheorem{cor}[theorem]{Corollary}
\newtheorem{lemma}[theorem]{Lemma}
\title{\bf A Dehn surgery description of regular finite cyclic covering spaces of rational homology spheres}

\author{ Cynthia L. Curtis }

\date{}

\begin{document}

\maketitle

\begin{abstract}
We provide related Dehn surgery descriptions for
rational homology spheres 
and a class of their regular finite cyclic covering spaces.  As an application,
we use  the surgery 
descriptions to relate the Casson invariants of the covering spaces to that
of the base space.  Finally, we  
show that this places restrictions on the number of finite and cyclic 
Dehn fillings of the knot complements in the covering spaces beyond those 
imposed by Culler-Gordon-Luecke-Shalen and Boyer-Zhang.
\end{abstract}


\noindent {\em AMS classification:} 57

\section{Introduction}

It is well-known that any closed, oriented, connected 
3-manifold may be obtained by Dehn surgery on a link
in  $S^3$.  In recent years, 3-manifold theorists have exploited this fact, 
computing certain 3-manifold invariants by describing how the invariant 
changes under Dehn surgery.   This has been an effective method for computing
invariants of individual manifolds.  However, until now, there was no known 
procedure for relating Dehn surgery descriptions of manifolds with those of
their covering spaces.  Therefore, although invariants for a manifold and a
covering space of the manifold could each be computed using Dehn surgery 
formulas, no general statements could easily be made regarding the 
relationship between the invariants. 

Let $(X, \tilde{X})$ be a 3-manifold pair, where $X$ is a rational 
homology sphere and 
$\tilde{X}$ is a regular finite cyclic covering space of $X$, say with
$\pi_1(X)/\phi_*(\pi_1(\tilde{X})) = {\displaystyle\Bbb Z} / 
k \displaystyle\Bbb Z$, where 
\mbox{$\phi:\tilde{X}\rightarrow X$} is the projection map.  Since 
this group is abelian, we may 
factor the quotient map 
$\pi_1(X)\rightarrow \pi_1(X)/\phi_*(\pi_1(\tilde{X}))$ through
the first homology group ${\rm H}_1(X;\displaystyle\Bbb Z)$.  

\begin{defn}

We call the covering $\tilde{X}\rightarrow X$ 
\mbox{\bf torsion-split} 
if there exists a homology decomposition
${\rm H}_1(X;{\displaystyle\Bbb Z}) = {\displaystyle\Bbb Z} 
/ kp {\displaystyle\Bbb Z} \oplus H$ (where possibly $H = 0$) 
satisfying:
\begin{itemize}
\item[i)] the decomposition is a decomposition of the torsion linking pairing on
${\rm H}_1(X)$: i.e. if $\alpha$ is a generator of the ${\displaystyle\Bbb Z} 
/ kp {\displaystyle\Bbb Z}$-summand and $\beta_1, \beta_2, \ldots, \beta_j$ is 
a torsion basis for $H$, then $link(\alpha,\alpha) = m/n$ for some $m$ 
relatively prime to $n$ and $link(\alpha,\beta_i) = 0$ for all $i$
\item[ii)] any generator of the ${\displaystyle\Bbb Z} 
/ kp {\displaystyle\Bbb Z}$-summand maps to a generator  
in $\pi_1(X)/\phi_*(\pi_1(\tilde{X}))$
under the quotient map 
$ {\rm H}_1(X) \rightarrow \pi_1(X)/\phi_*(\pi_1(\tilde{X}))$
\item[iii)] $H$ maps to 0 in $\pi_1(X)/\phi_*(\pi_1(\tilde{X}))$.
\end{itemize}

\end{defn}

In this paper, we provide a Dehn surgery description for torsion-split regular
$k$-fold cyclic covering space pairs 
$(X, \tilde{X})$ with base space a rational homology sphere.  
Specifically, we prove the following:

Let $K$ be a knot in $S^3$, and let $L =  (L_1,L_2,\ldots,L_n)$ be a link in 
$S^3$.  Assume $K$ bounds a Seifert surface $\Sigma$ with the property that 
there exist  Seifert surfaces 
$\Sigma_1, \Sigma_2,\ldots,\Sigma_n$ for $L_1,L_2,\ldots,L_n$ disjoint 
from a neighborhood of $\Sigma$.  Let $p$, $q$, and $k$ 
be integers with 
$1 \leq q$, $1 \leq |p|$, and $k > 1$.  Assume $kp$ and $q$ are relatively 
prime.  Let $M$ be the 3-manifold obtained by $kp/q$-Dehn surgery on $K$ in 
$S^3$ followed by surgery on $L$ 
with surgery coefficients $I = (i_1,i_2,\ldots,i_n)$, where $i_j = \pm 1$ for 
each $j$.

Note that for each  $j = 1,2,\ldots,n$, the component $L_j$ of $L$ has $k$ 
disjoint lifts in the $k$-fold branched cyclic cover of $S^3$ branched along 
$K$, since $(K,L_j)$ is a boundary link.
Let $\tilde{M}$ be the 3-manifold obtained by $p/q$-Dehn surgery on the lift 
of $K$ to the 
$k$-fold branched cyclic cover of $S^3$ branched along $K$, followed by 
surgery on the 
$k$ lifts of $L$ with surgery coefficient $i_j$ for every lift of $L_j$.  

\begin{theorem}
$\tilde{M}$ is a regular $k$-fold cyclic covering space of $M$. 
\end{theorem}

Call $(K,L,k,p,q,I)$ a {\em pairwise Dehn surgery description} for  
$(M,\tilde{M})$.  Then we also have

\begin{theorem}
Let $(X,\tilde{X})$ be a torsion-split regular $k$-fold cyclic covering space pair with base space $X$ a rational homology sphere.  Then 
$(X,\tilde{X})$ has a pairwise Dehn surgery description.
\end{theorem}

Thus, every torsion-split regular $k$-fold cyclic covering space pair over a rational homology sphere has a pairwise Dehn 
surgery description.  This, then, may be used for computing 3-manifold 
invariants for the pair.  
This should be a necessary first step in drawing general conclusions about 
the relationships between invariants for $X$ and those for $\tilde{X}$.

The paper is outlined as follows.  In Section 2, we show that the manifolds 
generated by the pairwise Dehn surgery description $(K,L,k,p,q,I)$ are a 
regular $k$-fold cyclic 
covering space pair.  In Section 3, we show that all torsion-split regular $k$-fold cyclic 
covering space pairs over rational homology spheres 
arise in this way.  
In Section 4, as an application,
we compute Casson invariants for the pair $(X,\tilde{X})$.  

In Section 5,
we discuss the implications of our work for generating 
examples of finite and cyclic Dehn fillings of 3-manifolds with toral boundary.
Specifically, let $K$ be a knot homologous to 0 in a closed, 
connected, oriented 3-manifold $Y$, and let $X$ be the result of $p/q$-Dehn 
surgery on $K$ in $Y$. We show that if $\pi_1(X)$ is finite (resp. cyclic), 
then so too is the fundamental group of the 3-manifold obtained by $r/q$-Dehn
surgery on the lift of $K$ to the $p/r$-fold cyclic branched cover of $Y$ 
branched along $K$, where $r$ is any positive integer dividing $p$.  Thus, 
one example of a finite or cyclic Dehn filling may generate many more such 
examples.  Moreover, this restricts the possible number of finite and cyclic 
surgery slopes in the cyclic covering spaces of $Y - nbhd(K)$.

Finally, in the appendix, we demonstrate some relationships between certain
Alexander polynomials which are required in Section 4.

We remark that the construction presented here may generalize to pairs 
$(X,\tilde{X})$ for which $X$ is not a rational homology sphere.  This theory,
together with applications, will be the topic of a future paper.

I am indebted to Andrew Clifford for his help throughtout the writing process
and to Steve Boyer for his corrections to an early version of this paper. 
I also thank Nancy Hingston and Andy Nicas
for their helpful comments.

\section{Generating covering space pairs}

Retain all notation from Section 1.  Further,
denote by $\overline{S}^3_K$ the $k$-fold cyclic covering space 
of $S^3$ branched over $K$.  
Let $N$ be the closed 3-manifold resulting from $kp/q$-Dehn surgery on $K$ in 
$S^3$, and let $\tilde{N}$ denote the 3-manifold which is the result of 
$p/q$-Dehn surgery on the lift of $K$ to $\overline{S}^3_K$.  Finally, if 
$Y$ is any 3-manifold and $K$ is a knot in $Y$, we abuse notation and denote 
by $Y-K$ the compact 3-manifold formed by removing an open tubular neighborhood
of $K$.

Before proving  Theorem 1.2, we prove the following:

\begin{lemma}
$\tilde{N}$ is a regular $k$-fold cyclic covering space of $N$.
\end{lemma}

\noindent
{\em Proof of Lemma 2.1:}  
Let $\overline{K}$ denote the lift of $K$ to $\overline{S}^3_K$, and let 
$\overline{m}$ denote its meridian.  Let $m$ denote the meridian of $K$.  Let 
$l$ denote the preferred longitude of $K$ on the boundary of $S^3 - K$, and 
let $\overline{l}$ be a longitude in the boundary of 
$\overline{S}^3_K - \overline{K}$ which 
lifts $l$.
Then $\tilde{N}$ is obtained by gluing a solid torus to 
$\overline{S}^3 - \overline{K}$, 
sending a meridian of the solid torus to the simple closed curve 
$\overline{m}^p \overline{l}^q$.  Similarly, $N$ is obtained by gluing a 
solid torus to $S^3 - K$, sending 
a meridian of the solid torus to $m^{kp}l^q$.

Let $\phi$ denote the regular $k$-fold cyclic covering map of the knot
complements.  We wish to show that $\phi$ extends to a regular $k$-fold cyclic
(unbranched) covering map  
$\tilde{N} \rightarrow N$.

Note that $\phi$ restricts to a regular $k$-fold cyclic covering map from the
boundary of $\overline{S}^3 - \overline{K}$ to $S^3 - K$.  Thus, the map from 
the boundary of the solid torus in $\tilde{N}$ to the boundary of the solid 
torus in $N$ is a regular $k$-fold cyclic covering map.

Now a $k$-fold cyclic covering
map from the boundary of one solid torus to the boundary of another solid
torus extends to a regular $k$-fold cyclic covering map
from the solid torus to the solid torus if and only if the meridian of the 
solid torus in the initial solid torus is taken to the meridian of the 
final solid torus.  In our case, the meridian of the solid torus in $\tilde{N}$
is $\overline{m}^p\overline{l}^q$, and the meridian of the solid torus in $N$
is $m^{kp}l^q$.  But $\phi$ takes $\overline{m}$ to $m^k$ and $\overline{l}$ to
$l$.  Hence, $\phi(\overline{m}^p\overline{l}^q) = m^{kp}l^q$.
It follows that $\phi$ takes a meridian of the solid torus in $\tilde{N}$ to 
a meridian of the solid torus in $N$, and therefore $\phi$ extends to a 
covering map $\tilde{N}\rightarrow N$.
$\Box$

We now prove Theorem 1.2.

\noindent {\em Proof of Theorem 1.2:}  We know that $\tilde{N}$ is a regular 
$k$-fold cyclic covering space of $N$ by Lemma 2.1. 
Let $\phi$ denote the covering map.
We abuse notation and denote the 
image of $L_j$ in $N$ by $L_j$ for $j = 1,2,\ldots,n$. Let 
$\overline{L}_j$ denote the inverse image under $\phi$ of $L_j$ for each $j$.  
Since $(K,L_j)$ is a boundary link in $S^3$ for each $j$, we see 
that $\overline{L}_j$ consists of $k$ disjoint simple closed curves in $\tilde{N}$.  
Choose a
Seifert surface $\Sigma$ for $K$ and Seifert surfaces $\Sigma_j$ for $L_j$ 
in $S^3$ such that $\Sigma \bigcap \Sigma_j = \emptyset$.  
Then the image of $\Sigma_j$ in 
$N$ is a Seifert surface for $L_j$ in $N$, and it lifts to $k$ Seifert surfaces
$\overline{\Sigma}_{j,1},\overline{\Sigma}_{j,2}, \ldots, 
\overline{\Sigma}_{j,k}$ for 
the $k$ lifts $\overline{L}_{j,1},\overline{L}_{j,2},\ldots,
\overline{L}_{j,k}$ of $L_j$ 
in $\overline{L}_j$.

Recall that $\overline{S}^3_K$ may be explicitly constructed according to the
following outline:  
let $Y^0 = \Sigma \times (-1,1)$ be an open bicollar of $\Sigma$, and let 
$Y = Y^0/ (K\times (-1,1) \sim K)$.  Let $Y^-$ be the manifold obtained 
by removing $K$ from $Y$.  Glue $k$ 
copies of $S^3 - \Sigma$ together along $k$ copies of $Y^-$, alternating copies
of $S^3 - \Sigma$ with copies of $Y^-$.  Finally, glue $K$ back in to compactify.
For a precise description of this
construction, see [R, pp 128 - 131 and pp 297-298].

Now since $(K,L_j)$ is a boundary link in $S^3$, we see that the $k$ lifts  of
$\Sigma_j$ to $\overline{S}^3_K$ are contained in the $k$ disjoint copies of 
$S^3 - Y$.  It follows that these lifts of $\Sigma_j$ are disjoint.  
Then clearly 
the Seifert surfaces $\overline{\Sigma}_{j,1}, \overline{\Sigma}_{j,2},
\ldots,\overline{\Sigma}_{j,k}$ in $\tilde{N}$ are also disjoint.

Now the Dehn surgery on $L_j$ may be 
carried out in a neighborhood $Z_j$ of $\Sigma_j$.  Moreover, 
if we take $Z_j$ to be 
sufficiently small, then the inverse image $\phi^{-1}(Z_j)$
consists of $k$ disjoint copies $\overline{Z}_{j,i}$ of $Z_j$ which are 
neighborhoods of the $\overline{\Sigma}_{j,i}$. Clearly the Dehn surgery 
on $L_j$ in 
$Z_j$ induces an identical  Dehn
surgery on $\overline{L}_{j,i}$ in $\overline{Z}_{j,i}.$
Therefore every point of  $M$ has $k$ distinct inverse images in
$\tilde{M}$, each with a neighborhood which is carried homeomorphically to a 
neighborhood of the point in $M$.  Thus, $\tilde{M}$ is a $k$-fold covering 
space of $M$.

Finally, note that for each $j$, the automorphism group 
${\displaystyle\Bbb Z}/ k {\displaystyle\Bbb Z}$ of 
$\phi:\tilde{N} \rightarrow N$ cyclically permutes the $k$ disjoint copies of the 
Seifert surface for 
$L_j$ in $\tilde{N}$.  Clearly, then, the automorphism group of the 
covering space 
$\tilde{M}\rightarrow M$ is also ${\displaystyle\Bbb Z} / k {\displaystyle\Bbb Z}$,
and the images of the $k$ lifts of 
$\overline{L}_{j,i}$ after Dehn surgery are permuted by the 
automorphism group.  Since 
the covering space $\tilde{N}\rightarrow N$ was regular, it follows that 
$\tilde{M}\rightarrow M$ is regular.
$\Box$

We remark that the covering $\tilde{M} \rightarrow M$ is clearly torsion-split by construction.

\section{Completeness of the construction}

In this section, we show that the construction described in Section 
1 is complete in that it 
generates all torsion-split regular $k$-fold cyclic covering space pairs 
$(X,\tilde{X})$ over rational homology spheres. We prove:

\noindent {\bf Theorem 1.3} {\em Let $(X,\tilde{X})$ be a torsion-split
regular $k$-fold 
cyclic covering space pair with base space $X$ a rational homology sphere.  Then $(X,\tilde{X})$ has a pairwise 
Dehn surgery description.}

We first prove the following lemma, which is a straightforward 
generalization of a lemma 
of S. Boyer and D. Lines [BL].

\begin{lemma}

Let $W$ be a rational homology sphere 
with ${\rm H}_1(W;{\displaystyle\Bbb Z}) = 
{\displaystyle\Bbb Z} / n {\displaystyle\Bbb Z} \oplus H$ for some 
finite abelian group $H$.  Assume further that the homology decomposition 
arises as a decomposition of the torsion linking pairing on ${\rm H}_1(W)$.
Then there is a 3-manifold 
$V$ with ${\rm H}_1(V;{\displaystyle\Bbb Z}) = H$, a knot ${\cal K}$ 
homologous to 0 in $V$, and an 
integer $m$ such that $W$ is the result of $n/m$-Dehn surgery on ${\cal K}$ 
in $V$.
\end{lemma}

\noindent {\em Proof of Lemma 3.1:}
Let $\alpha$ be a generator of the ${\displaystyle\Bbb Z} / n 
{\displaystyle\Bbb Z}$ summand of 
${\rm H}_1(W;\displaystyle\Bbb Z)$.  
Represent $\alpha$ by a curve $C$ in $W$.

Note that the torsion subgroup of  
${\rm H}_1(W - C;\displaystyle\Bbb Z)$ is just $H$.  
To see this, note that by assumption, $link(\alpha,\alpha)= t/n$ for some 
integer $t$ relatively
prime to $n$, where $link(\mbox{\_},\mbox{\_})$ is 
the torsion linking pairing on ${\rm H}_1(W)$.
Therefore $C$ intersects the surface with boundary $nC$, and 
$nC$ does not bound in $W - C$.  On the other hand, there exist 
generators $\beta_i$ of $H$ such that 
$link(\alpha,\beta_i) = 0$.  
It follows that $\beta_i$ has the same finite order
in ${\rm H}(W - C)$ as in ${\rm H}(W)$.

Let $T(C)$ be a tubular neighborhood of $C$.  Then 
${\rm H}_2(W ,W - T(C);\displaystyle\Bbb Z)$ is infinite 
cyclic with generator a meridional disk of $T(C)$.  From the exact sequence
\[0\rightarrow {\rm H}_2(W ,W - T(C);{\displaystyle\Bbb Z}) \rightarrow 
{\rm H}_1(W - T(C);{\displaystyle\Bbb Z}) \rightarrow 
{\rm H}_1(W;{\displaystyle\Bbb Z}) \rightarrow 0\]
it follows that ${\rm H}_1(W - T(C);{\displaystyle\Bbb Z}) \cong 
{\displaystyle\Bbb Z} \oplus H$.

Let $C'$ be a simple closed curve on $\partial T(C)$ generating 
the infinite cyclic 
summand of ${\rm H}_1(W - T(C);\displaystyle\Bbb Z)$.  
Attach a solid torus to $W - T(C)$ sending the meridian to $C'$; 
call the resulting manifold $V$.  Let ${\cal K}$
denote the core of the surgery torus.  Then 
${\rm H}_1(V;{\displaystyle\Bbb Z}) \cong  H$, and 
$W$ is the result of $n/m$-Dehn surgery on ${\cal K}$ for some integer $m$.  
Moreover, ${\cal K}$ is homologous to 0 in $V$, 
since all curves on $\partial T(C)$ represent 0 in $H$.
$\Box$

\noindent We now prove the theorem.

\noindent {\em Proof of Theorem 1.3:}
Let $\tilde{X}\rightarrow X$ be a torsion-split regular $k$-fold cyclic covering
space with $X$ a rational homology sphere.  Let ${\rm H}_1(X;{\displaystyle\Bbb Z}) = {\displaystyle\Bbb Z} 
/ kp {\displaystyle\Bbb Z} \oplus H$ be a homology decomposition for $X$ obeying
properties {\em i) - iii)} of Definition 1.1.
We may apply Lemma 3.1 
to find a 3-manifold $V$ with ${\rm H}_1(V;{\displaystyle\Bbb Z}) = H$, 
a knot ${\cal K}$ homologous to 0 in 
$V$, and an integer $q$ such that $X$ is the result of 
$kp/q$-Dehn surgery on ${\cal K}$ in 
$V$. Let $\Sigma$ be a Seifert surface for ${\cal K}$ in $V$.

It is well-known that there exists a link   
${\cal L}=({\cal L}_1,{\cal L}_2,\ldots,{\cal L}_n)$ in $V$ such that $S^3$ 
is the result of surgery on ${\cal L}$.  Let 
$L = (L_1,L_2,\ldots,L_n)$ be the image of 
${\cal L}$ in $S^3$.  Then $V$ may be 
obtained from $S^3$ by surgeries on the components of $L$. Moreover, we may choose 
${\cal L}$ in such a way that the surgery coefficients for the components
of $L$ are all $\pm 1$.

Choose $\alpha_1,\alpha_2,\ldots,\alpha_{2g}$ a collection of simple closed curves
on $\Sigma$ representing a homology basis for
$\Sigma$.  Isotope ${\cal L}$ without changing any crossings of ${\cal L}$
so that the linking number of ${\cal L}_i$ and $\alpha_j$ is 0 for any
$i = 1,2,\ldots,n$ and any $j = 1,2,\ldots,2g$ and so that ${\cal L}_i \cap \Sigma = \emptyset$.  
We continue to denote by
${\cal L}$ and $L$ the images of ${\cal L}$ and $L$ under the isotopy.
Let $K$ be the image of ${\cal K}$ in $S^3$ after surgery on ${\cal L}$.
Abusing notation yet again, we denote by $\Sigma$ and $\alpha_j$ the image of $\Sigma$ and $\alpha_j$ in
$S^3$.  Note that the linking number of $L_i$ and $\alpha_j$ is 0 for any
$i = 1,2,\ldots,n$ and any $j = 1,2,\ldots,2g$ by Lemma A.5. Furthermore $L_i$ does not meet $\Sigma$.

We require the following

\begin{lemma}
Let $C$ be a knot in $S^3$, and let $S$ be a Seifert surface for
$C$.  Let $x_1,x_2,\ldots,x_{2g}$ be a collection of curves on $S$ representing a
basis for ${\rm H}_1(S;{\displaystyle\Bbb Z})$.  Let $D$ be a knot in $S^3$ such that
$D \cap S = \emptyset$ and 
$lk(x_i,D)=0$ for $i = 1,2,\ldots,2g$.  Then $D$ bounds a Seifert surface
disjoint from $S$.
\end{lemma}

\noindent {\em Proof of lemma:}  Let $S'$ be a Seifert surface
for $D$ meeting $S$ transversely.  If $S' \cap S = \emptyset$, we are done.  Otherwise, 
$S' \cap S$ is a collection $y_1,y_2,\ldots,y_m$ of oriented simple closed curves.  Now
the homology class represented by $y_1+y_2+\ldots +y_m$ in ${\rm H}_1(S;{\displaystyle\Bbb Z})$ 
must be 0.  For if this
class were non-zero, then there would be a curve $x_i$ on $S$ whose oriented
intersection number with $y_1+y_2+\ldots +y_m$ was non-zero.  Then
$x_i$ would have non-zero intersection number with $S'$ and hence
would have non-zero linking number with $D$.    

Now remove  from
$S'$ the components of $S' - S$ which do not meet $D$. 
The resulting surface $S''$ has boundary $y_1 \cup y_2 \cup \ldots \cup y_m \cup D$. 
Since the sum of the $y_i$'s represents 0 in ${\rm H}_1(S;{\displaystyle\Bbb Z})$, 
the curves $y_i$ cobound a collection
of subsurfaces of $S$.  Then we may glue a collection of parallel copies
of these surfaces to the appropriate boundary components of $S''$ to form
a new two-sided surface $S'''$ with boundary $D$. Pushing the parallel
copies of subsurfaces of $S$ apart and away from $S$, we will find that $S'''$ is
embedded and disjoint from $S$. $\Box$

Now returning to the proof of the theorem:  

Applying the lemma to each of the link components $L_i$, we may find Seifert
surfaces $\Sigma_1,\Sigma_2,\ldots,\Sigma_n$ for $L_1,L_2,\ldots,L_n$,
respectively, such that $\Sigma \cap \Sigma_i = \emptyset$ for 
$i =1,2,\ldots,n$.  We show

\noindent {\em Claim:}  $(K,L,k,p,q,I)$ is a pairwise 
Dehn surgery description for $(X,\tilde{X})$.

\noindent {\em Proof of Claim:}  It is clear from the construction 
that $(K,L,k,p,q,I)$ is a 
pairwise Dehn surgery description for some pair $(X,X')$ with base space $X$.
  Let $\psi: X' \rightarrow X$ be the projection map.  Then 
$\psi_*(\pi_1(X'))$ is the kernel of the 
homomorphism 
\[\pi_1(X)\rightarrow {\rm H}_1(X;{\displaystyle\Bbb Z}) 
\rightarrow {\displaystyle\Bbb Z} / k {\displaystyle\Bbb Z},\] 
where the latter homomorphism is the projection 
\[{\rm H}_1(X;{\displaystyle\Bbb Z}) \rightarrow 
{\rm H}_1(X;{\displaystyle\Bbb Z}) / H = 
{\displaystyle\Bbb Z} / kp {\displaystyle\Bbb  Z} \rightarrow 
{\displaystyle\Bbb Z} / k {\displaystyle\Bbb Z}.\]  But this 
kernel is precisely $\phi_*(\pi_1(\tilde{X}))$.  Hence $X' = \tilde{X}$.
$\Box$

\section{  Casson-Walker invariants for pairs}
In 1985, Andrew Casson defined an invariant $\lambda$ for integral homology 
3-spheres.  Roughly, this invariant counts the signed equivalence classes 
of SU(2)-representations of the fundamental group of the 3-manifold.  
This invariant was extended 
to an invariant for oriented rational homology 3-spheres by Kevin Walker 
in [W], and Christine 
Lescop derived a combinatorial formula extending the invariant to 
arbitrary closed, oriented 3-manifolds 
in [L].

A number of mathematicians have explored the Casson-Walker 
invariant for branched 
covers of links in $S^3$.  David Mullins computed the invariant for 2-fold 
branched covers in the case when the 2-fold branched cover is a 
rational homology sphere in [M].  More general results for $k >2$ may be found in [GR].
The invariants for $n$-fold branched covers $\overline{S}^3_K$
of particular families of knots have been 
computed by J. Hoste, A. Davidow, and K. Ishibe in [H], [D], and [I]. 
Garoufalidis generalizes several of these formulas in [G].

We study the Casson-Walker invariants of pairs $(X,\tilde{X})$.  
As in the early sections of the paper, 
we assume $\tilde{X} \rightarrow X$ is a torsion-split regular $k$-fold cyclic 
covering over a rational homology sphere.  In this section, we assume further that
$\tilde{X}$ is a rational homology sphere and that $(X,\tilde X)$ has a Dehn surgery
description $(K,L,k,p,q,I)$ with $\overline S^3_K$ a rational homology
sphere.

In what follows, for any pair of non-zero integers
$x$ and $y$ which are relatively prime, let 
$s(x,y)$ denote the Dedekind sum defined by
\[s(x,y) = sign (y)\sum_{j = 1}^{|x|}((j/y))((jx/y))\]
where
\[((z)) = \left\{ \begin{array}{ll}
	0 & z \in {\displaystyle\Bbb Z}\\
	z - [z] - 1/2 & \mbox{else}
	\end{array}
\right. \]
For a knot $C$ in a rational homology sphere, let $\Delta_C$ denote the Alexander
polynomial of C, normalized so that it is symmetric in $t^{1/2}$ and $t^{-1/2}$ and so that $\Delta_C(1)= 1$.
We show

\begin{theorem}
Let $(K,L,k,p,q,I)$ be a pairwise Dehn surgery 
description for $(X,\tilde{X})$.  Then
\[\lambda(\tilde{X}) = 
k \lambda(X)  + q/p (\Delta^{''}_{\overline{K}}(1) - \Delta_K^{''}(1)) - k s(q,kp) + s(q,p) + \lambda(\overline{S}^3_K).\]
Here, $\overline{S}^3_K$ denotes 
the $k$-fold branched cyclic cover of $S^3$ branched along $K$ and $\overline{K}$ 
denotes the lift of $K$ to $\overline{S}^3_K$, as above.
\end{theorem}

Note that $\Delta_{\overline{K}}$ can be computed from $\Delta_K$.  This relationship is described in the appendix.
We remark further that in the case $k=2$, the invariant $\lambda(\overline{S}^3_K)$ 
can be computed whenever $\tilde{X}$ is a rational homology sphere 
using the work of 
Mullins.  For $k > 2$, the results of Garoufalidis and Rozansky apply.  For certain families of  knots,  the invariant 
$\lambda(\overline{S}^3_K)$ can be computed 
for any value of $k$ using the results of Hoste, Davidow, Ishibe, and Garoufalidis.

{\em Proof of Theorem 4.1:} Retain all notation from previous sections.  
We begin by noting that 
\begin{equation}
\lambda(N) = (q/kp) \Delta^{''}_K (1) + s(q,kp)
\end{equation}
and
\begin{equation}
\lambda(\tilde{N}) = 
\lambda(\overline{S}^3_K) + (q/p) \Delta^{''}_{\overline{K}}(1) + s(q,p)
\end{equation}
by Proposition 6.2 of [W], since $K$ is a knot in $S^3$ and 
since $\overline{K}$ is homologous to 0 in $\overline{S}^3_K$.

Now $X$ is obtained from $N$ by $I$-surgery on $L$.  We know that $\lambda(X) 
- \lambda(N)$ may be obtained using the surgery formulae developed by Walker 
in [W].  These formulae depend on the coefficients $I$, as well as the link
$L$ and the Alexander polynomials of the components of L.  Similarly $\tilde{X}$ is obtained from $\tilde{N}$ by $I$-surgery on 
each lift $\overline{L}_j$ of $L$ in $\tilde{N}$.  

Now by Proposition A.7 in the appendix, the Alexander polyomial of each component of  $\overline{L}$
is equal to that of the corresponding component of $L$.  Therefore, since a neighborhood of 
each component $\overline{L}_{j,i}$ of $\overline{L}_j$ in $\tilde{N}$ 
is carried homeomorphically onto a 
neighborhood of $L_j$ in $N$ and the coefficients $I$ correspond, it is clear that for any $j$, the Casson-Walker
 invariant of the manifold $\tilde{N}_j$ obtained by doing $I$-surgery on 
$\overline{L}_j$ in $\tilde{N}$ obeys the formula $\lambda(\tilde{N}_j) - 
\lambda(\tilde{N}) = \lambda(X) - \lambda(N)$.  But the lifts 
$\overline{\Sigma}_{j,i}$
of $\Sigma_j$ are contained in disjoint copies of $S^3 - Y$ in $\tilde{N}$,
as shown in the proof of Theorem 1.2,   
so $\overline{\Sigma}_{j,i} \cap \overline{\Sigma}_{l,m} = \emptyset$ if 
$i \neq m$.  Then also $\lambda(\tilde{N_{jk}}) - \lambda(\tilde N_j) = 
\lambda(\tilde N_j)-\lambda(N)$, where $\tilde N_{jk}$ is the manifold obtained
by doing $I$-surgery on the image of $\overline{L}_k$ in $\tilde N_j$.
It follows that 
\begin{equation}
\lambda(\tilde{X})-\lambda(\tilde{N}) = k
(\lambda(X) - \lambda(N)).
\end{equation} 
 The theorem follows after suitably combining equations (1) - (3).
$\Box$

We remark that analogous results for the generalizations of $\lambda$ counting 
representations in SO(3), U(2), Spin(4), and SO(4) may be immediately obtained 
using the results of [C].

\section{Cyclic and finite Dehn surgeries}
In recent years, many new and exciting results have come to light 
concerning Dehn fillings of 3-manifolds with toral boundary.  Among these are
the Cyclic Surgery Theorem of Culler, Gordon, Luecke, and Shalen [CGLS] and 
the work of Boyer and Zhang concerning finite fillings [BZ1] and [BZ2].  

We review some basic definitions. 

\begin{defn}
Let $K$ be a knot homologous to 0 in an oriented 3-manifold $Y$. A \mbox{\bf
slope} 
of $K$ is the unoriented isotopy class of a non-trivial simple closed curve in 
$\partial(Y-K)$.  The \mbox{\bf distance} between two slopes is their 
geometric intersection number.
\end{defn}

Recall that for any knot homologous to 0 in an oriented 3-manifold, the set 
of slopes of the knot is canonically isomorphic to 
${\overline{\displaystyle\Bbb Q}}$.  Thus, we denote by $p/q$ the slope which 
is the isotopy class of a curve which is $p$ times a meridian plus $q$ times a 
longitude.

\begin{defn}
With $K$ and $Y$ as above, suppose that
$p/q$ is a slope of $K$ such that the manifold $X$ which is the result 
of $p/q$ Dehn surgery on $K$ in $Y$ has finite fundamental group.  
Then we call $p/q$ a
\mbox{\bf finite surgery slope} of $K$.  Similarly, if $X$ has cyclic 
fundamental
group, we call $p/q$ a \mbox{\bf cyclic surgery slope} of $K$.
\end{defn}

The results of Culler, Gordon, Luecke, and Shalen and of Boyer and Zhang 
state that for most irreducible knot complements $Y-K$ as above, 
there are at most 
3 cyclic surgery slopes of $K$ and at most
6 finite surgery slopes.  Moreover
these slopes may be distance at most 1 apart in the cyclic case 
and at most 5 apart in the finite case.  A detailed survey of work in this area may be found in [B].

Here, we note that given a knot $K$ homologous to 0 in a closed, oriented
3-manifold $Y$ and a finite or cyclic surgery slope of 
$K$, our work leads to explicit
examples of more such fillings. 
Specifically, let $\overline{Y}_K$ 
denote the $k$-fold branched cyclic cover $Y$ branched along $K$, and
let $\overline{K}$ denote the lift of $K$ to $\overline{Y}_K$.  We show

\begin{theorem}
Let $K$ be a knot homologous to 0 in a closed oriented 3-manifold $Y$. 
Let $k$, $p$,
and $q$ be integers with $k > 1$, $|p| \geq 1$, and 
$q \geq 1$.  Then 
$kp/q$ is a finite surgery slope of $K$ if and only if $p/q$ is a 
finite surgery 
slope of $\overline{K}$ in $\overline{Y}_K$.
Moreover if $kp/q$
is a cyclic surgery slope of $K$, then $p/q$ is a 
cyclic surgery slope of
$\overline{K}$.  Finally, if $p/q$ is a 
cyclic surgery slope of
$\overline{K}$ and $p \neq 1$, then $kp/q$
is a cyclic surgery slope of $K$.
\end{theorem}

\noindent {\em Proof:}  Let $X$ denote the manifold resulting from $kp/q$-Dehn
 surgery on $K$, and let $\tilde{X}$ denote the result of $p/q$-Dehn surgery
on $\overline{K}$.  As in the proof of Theorem 1.3, we
may find a knot $C$ and a link $L = (L_1,L_2,\ldots,L_n)$ in $S^3$ satisfying
\begin{itemize}
	\item there exist Seifert surfaces $\Sigma$ and $\Sigma_1,
\Sigma_2,\ldots,\Sigma_n$ for $C$ and $L_1,L_2,\ldots,L_n$, respectively, with
$\Sigma\cap\Sigma_j = \emptyset$ for $j=1,2,\ldots,n$ 
	\item $Y$ is the result of $I$-surgery on $L$ 
 for some $I=(i_1,i_2,\ldots,i_n)$ with $i_j = \pm 1$
	\item $K$ is the image of $C$ in $Y$.
\end{itemize}
(Note that none of these steps requires $Y$ to be a rational homology sphere.)

Now $(C,L,k,p,q,I)$ is a Dehn surgery description for some pair $(Z,\tilde Z)$.
But since $(C,L_i)$ is a boundary link for $i=1,2,\ldots,n$, it is clear that 
$kp/q$-surgery on $C$ followed by $I$-surgery on the image of $L$ yields the same manifold
as $I$-surgery on $L$ followed by $kp/q$-surgery on the image of $C$. Hence $X \cong Z$.
Also, since $(C,L_i)$ is a boundary link, we see that $\overline Y_K$ may be obtained by $I$-surgery
on the $k$ lifts of $L$ to $\overline S^3_C$.  Moreover $p/q$-surgery on $\overline C$
in $\overline S^3_C$ followed by $I$-surgery on the images of each of the $k$ lifts of
$L$ to $\overline S^3_C$ yields the same manifold as $I$-surgery on each of the lifts of
$L$ to $\overline S^3_C$ followed by $p/q$-surgery on the image $\overline K$ of $\overline C$
in $\overline Y_K$.  Hence $\tilde X \cong \tilde Z$.

Now it follows from Theorem 1.2 that
 $\tilde{X}$ is a regular $k$-fold cyclic covering space of $X$.  Therefore the
fundamental group of
$\tilde{X}$ is an index $k$ subgroup of that of $X$.  
Then clearly $\pi_1(X)$ is finite if and only if $\pi_1(\tilde{X})$ is 
finite, and 
$\pi_1(\tilde{X})$ is cyclic if $\pi_1(X)$ is cyclic.

Finally, if $\pi_1(\tilde{X})$ is cyclic and $p \neq 1$, so 
$\pi_1(\tilde{X}) = {\displaystyle\Bbb Z}/p {\displaystyle\Bbb Z}$, 
then $\pi_1(X) = {\rm H}_1(X) =
 {\displaystyle\Bbb Z}/ k p {\displaystyle\Bbb Z}$, 
since ${\rm H}_1(X; {\displaystyle\Bbb Z})$ has a 
${\displaystyle\Bbb Z}/ k p {\displaystyle\Bbb Z}$ 
summand and $\pi_1(\tilde{X})$ is 
index $k$ in $\pi_1(X)$.  This proves the theorem.  
$\Box$

In light of [CGLS], [BZ1], and [BZ2], then, finding a finite or cyclic surgery slope 
$p/q$ for
some knot $K$ homologous to 0 in an oriented 3-manifold $Y$ severely restricts 
the set of possible finite and cyclic surgery slopes not only of $K$ itself, 
but also of lifts of $K$ to branched covers of $Y$ branched along $K$ of 
all orders dividing $p$.  For example, it is a result of Fintushel and Stern 
[FS] that 18- and 19-surgeries on the 
$(-2,3,7)$ pretzel knot yield lens spaces.  It follows that 1-surgery on the
lift of the $(-2,3,7)$ pretzel knot to either the 18- or 19-fold branched 
cover of $S^3$ branched along the pretzel knot yields $S^3$.  
Further, we find that 
$p$-surgery on the lift of the knot to the 
$(18/p)$-fold branched cover of $S^3$ 
branched along the knot also yields a lens space for $p = 2$, 3, 6, or 9.

Finally, we note that the results of [CGLS], [BZ1], and [BZ2] may be strengthened for 
manifolds of 
the form $\overline{Y}_K - \overline{K}$.  We have 
the following corollary, where 
$Y$, $K$, $\overline{Y}_K$, and $\overline{K}$ are defined as above:

\begin{cor}
Suppose $Y - K$ is irreducible and is not a Seifert fibered 
space.  Then $\overline{K}$ 
has at most one cyclic surgery slope $p/q$ with $p \neq 1$.

If also $Y - K$ is not a cable on the twisted $I$-bundle over the 
Klein bottle, then the 
distance between any two finite surgery slopes of $\overline{K}$ is 
at most $5/k$.

If $Y - K$ is also hyperbolic, then the distance between a cyclic 
surgery slope $p/q$ with $p \neq 1$ and any 
finite surgery slope is at most $2/k$.  Moreover the distance between any two 
finite surgery slopes is at most $3/k$. 
\end{cor}

\noindent{\em Proof:}  By Theorem 5.3, for any $p$ and $q$ with $|p|\geq 1$ 
and $q\geq 1$, 
we know that $p/q$ is a finite 
surgery slope of 
$\overline{K}$ if and only if 
$kp/q$ is a finite 
surgery slope of $K$ and that $kp/q$ is a cyclic surgery slope of
$K$ if $p/q$ is a cyclic surgery slope of $\overline{K}$ with $p \neq 1$. 
The distance between $kp'/q'$ and 
$kp/q$ is $k$ times the distance between $p'/q'$ and $p/q$.  
The assertions follow from 
the restrictions on the distance between cyclic (resp. finite) 
surgery slopes $kp/q$ for $K$ 
imposed by [CGLS], [BZ1], and [BZ2].
$\Box$

\appendix
\section{Alexander polynomials in $\overline S^3_K$}
In this section, we relate the Alexander polynomials $\Delta_{\overline K}$ and 
$\Delta_{\overline L_{j,i}}$ in $\overline S^3_K$ to $\Delta_K$ and
$\Delta_{L_j}$, respectively.  These results are used in Section 4.

We begin by recalling the definition of the Alexander polynomial of a knot 
in a rational homology sphere.  Details can be found in [M].

Let $C$ be a knot in a rational homology sphere, and denote by $Z$ the complement 
of a tubular neighborhood of $C$.  Further, denote by $\tilde Z$ the infinite
cyclic cover of $Z$ determined by 
$\pi_1(Z) \rightarrow {\rm H}_1(Z;{\displaystyle\Bbb Z})/torsion$.  Then
${\rm H}_1(\tilde Z;{\displaystyle\Bbb Q})$ is a module over 
${\displaystyle\Bbb Q}[{\displaystyle\Bbb Z}]$.

Now ${\displaystyle\Bbb Q}[{\displaystyle\Bbb Z}]$ is a principal ideal domain, so
\[{\rm H}_1(\tilde Z;{\displaystyle\Bbb Q}) \cong {\displaystyle\Bbb Q}[{\displaystyle\Bbb Z}]/(p_1) \oplus
{\displaystyle\Bbb Q}[{\displaystyle\Bbb Z}]/(p_2) \oplus
\ldots \oplus {\displaystyle\Bbb Q}[{\displaystyle\Bbb Z}]/(p_r). \]
We define the {\em order} of ${\rm H}_1(\tilde Z;{\displaystyle\Bbb Q})$ to
be the product ideal $(p_1p_2\ldots p_r)$.

\begin{defn}
An \mbox{\bf Alexander polynomial} $\Delta_C(t)$ of $C$ is any polynomial generating the order of 
${\rm H}_1(\tilde Z;{\displaystyle\Bbb Q})$.  We also call $\Delta_C(t)$ an 
\mbox{\bf Alexander polynomial} of the ${\displaystyle\Bbb Q}[{\displaystyle\Bbb Z}]$-module
${\rm H}_1(\tilde Z;{\displaystyle\Bbb Q})$.
\end{defn}

This is well-defined up to multiplication by polynomials of the form $c t^k$, where $c \in {\displaystyle\Bbb Q}$,
where $t$ generates the deck transformations of the covering $\tilde{Z} \rightarrow Z$, and where 
$k \in {\displaystyle \Bbb Z}$.  
Henceforth we write $p(t) \equiv q(t)$ if $p(t)$ and $q(t)$ are polynomials
with $p(t) = ct^kq(t)$ for some $c \in {\displaystyle\Bbb Q}$ and some 
$k \in {\displaystyle\Bbb Z}$.
 
We first relate $\Delta_{\overline K}$ to $\Delta_K$.  We prove the following theorem, which was pointed out to me by Steve Boyer.
The proof offered is that of Boyer.

\begin{prop} (Boyer) 
\[\Delta_{\overline K}(t^k) \equiv \prod_{j=0}^{k-1} \Delta_K (\zeta^j t)\]
where $\zeta$ is a primitive $k$th root of unity.
\end{prop}

\noindent {\em Proof (Boyer):}
Let $A$ be the matrix of multiplication by $t$ in ${\rm H}_1(\tilde Z;{\displaystyle\Bbb Q})$
with respect to the ${\displaystyle\Bbb Q}$-vector space structure.  We first show

\begin{lemma} (Boyer)
$\Delta_K(t) \equiv |A - tI|$
\end{lemma}

\noindent {\em Proof of lemma (Boyer):} 
Since ${\rm H}_1(\tilde Z;{\displaystyle\Bbb Q})$ decomposes as 
\[{\displaystyle\Bbb Q}[t,t^{-1}]/p_1(t) \oplus
{\displaystyle\Bbb Q}[t,t^{-1}]/p_2(t) \oplus
\ldots \oplus {\displaystyle\Bbb Q}[t,t^{-1}]/p_r(t),\]
we see that $A = \oplus_{j=1}^r A_j$, where $A_j$ is the matrix of multiplication 
by $t$ in ${\displaystyle\Bbb Q}[t,t^{-1}]/p_j(t)$.  Thus it suffices to show that
\[|A_j - tI| \equiv  p_j(t).\]  
But writing \mbox{$p_j(t) = b_0 + b_1 t + b_2 t^2 + \ldots + b_{s-1} t^{s-1} + t^s$},
we see that ${\displaystyle\Bbb Q}[t,t^{-1}]/p_j(t)$ has basis ${1,t,t^2,\ldots,t^{s-1}}$. 
Then 
\[ A_j = \left [ \begin{array}{lllll}
		0 & 0 & \ldots & 0 & -b_0\\
		1 & 0 & \ldots & 0 & -b_1\\
		0 & 1 & \ldots & 0 & -b_2\\
		\ldots & \ldots & \ldots & \ldots & \ldots\\
		0 & 0 & \ldots & 1 & -b_{s-1}
	\end{array} \right ] \]
It follows that $|A_j - tI| = (-1)^s p_j(t)$.
$\Box$

Now define a new ${\displaystyle\Bbb Q}[t,t^{-1}]$-module structure on 
${\rm H}_1(\tilde Z;{\displaystyle\Bbb Q})$ with mutliplication \mbox{$t\ast m = t^k m$}.
Clearly ${\rm H}_1(\tilde Z;{\displaystyle\Bbb Q})$ with the new multiplication is again a finitely generated, torsion
${\displaystyle\Bbb Q}[t,t^{-1}]$-module and hence a finite dimensional ${\displaystyle\Bbb Q}$-vector space.
We show

\begin{lemma}(Boyer) Let $\Delta_k(t)$ denote an Alexander polynomial for
${\rm H}_1(\tilde Z;{\displaystyle\Bbb Q})$ with the new
${\displaystyle\Bbb Q}[t,t^{-1}]$-module structure.  Then
\[ \Delta_k(t^k) \equiv \prod_{j=0}^{k-1} \Delta(\zeta^j t)\]
\end{lemma}

\noindent{\em Proof of lemma (Boyer):}
If $A$ is the matrix of multiplication by $t$ in 
${\rm H}_1(\tilde Z;{\displaystyle\Bbb Q})$ with the original 
${\displaystyle\Bbb Q}[t,t^{-1}]$-module structure,
then $A^k$ is the matrix of multiplication by $t$ in
${\rm H}_1(\tilde Z;{\displaystyle\Bbb Q})$ with the new 
${\displaystyle\Bbb Q}[t,t^{-1}]$-module structure.
Therefore by the previous lemma,
\begin{eqnarray}
\Delta_k(t^k) & \equiv  & |A^k - t^k I| \\
              & =  & \prod_{j=0}^{k-1} |A - \zeta^j t I|\\
              & =  & \prod_{j=0}^{k-1} \Delta(\zeta^j t)
\end{eqnarray}
This proves the lemma.
$\Box$

But now note that an Alexander polynomial $\Delta_{\overline K}(t)$ is an
Alexander polynomial of ${\rm H}_1(\tilde Z;{\displaystyle\Bbb Q})$ with the
second ${\displaystyle\Bbb Q}[t,t^{-1}]$-module structure, 
while $\Delta_K(t)$ is an Alexander polynomial
of ${\rm H}_1(\tilde Z;{\displaystyle\Bbb Q})$ with the original 
${\displaystyle\Bbb Q}[t,t^{-1}]$-module structure.  The theorem follows.
$\Box$

Thus, renaming $u = t^n$, say,  we obtain $\Delta_{\overline{K}}(u)$ 
in terms of $\Delta_K$.

This can be further simplified as follows: write 
\begin{equation} 
\Delta_K(t) = c_0 + c_1(t + t^{-1}) + c_2(t^2 + t^{-2}) + \ldots + c_n(t^n + t^{-n}).
\end{equation}
Applying Boyer's theorem to $(7)$ and noting that all cross-terms cancel, we see that
\begin{equation}
\Delta_{\overline K}(t^k) \equiv c_0^k + c_1^k \prod_{j=0}^{k-1}(\zeta^j t + \zeta^{-j} t^{-1}) + c_2^k \prod_{j=0}^{k-1}(\zeta^{2j} t^2 + \zeta^{-2j} t^{-2})
+ \ldots + c_n^k \prod_{j=0}^{k-1}(\zeta^{nj} t^n + \zeta^{-nj} t^{-n}).
\end{equation}
In particular, if $k$ is odd, note that $\prod_{j=0}^{k-1}(\zeta^{rj} t^r + \zeta^{-rj} t^{-r}) = t^{rk}+ t^{-rk}$ for any $r$.
Then equation (8) becomes
\begin{equation}
\Delta_{\overline K}(t^k) \equiv c_0^k + c_1^k(t^k+t^{-k}) + c_2^k(t^{2k} + t^{-2k}) + \ldots + c_n^k(t^{nk} + t^{-nk})
\end{equation}
and finally
\begin{equation}
\Delta_{\overline K}(u) \equiv c_0^k + c_1^k(u + u^{-1}) + c_2^k(u^2 + u^{-2}) + \ldots + c_n^k(u^n+u^{-n}).
\end{equation}

We now turn to relating $\Delta_{\overline L_{j,i}}$ to $\Delta_{L_j}$.
We return to the convention that $\Delta_C$ is symmetric in $t^{1/2}$ and $t^{-1/2}$ and that 
$\Delta_C(1)=1$, so that $\Delta_C$ is uniquely defined.

For the remainder of the paper, for any rational homology sphere $M$, 
let $lk_M(\_,\_)$ denote the linking number in $M$.
\begin{lemma}
Given knots $C$ and $C'$ in a rational homology sphere $M$  such that \mbox{$lk_M(C,C')=0$.}  
Let $\hat M$ denote the manifold resulting from $p/q$-Dehn surgery on $C'$ in $M$.  
Then 
\mbox{$lk_{\hat M}(\hat C,\hat D) = lk_M(C,D)$} for any knot $D$ in $M$, where $\hat C$ and $\hat D$ are the images of $C$ and
$D$ in $\hat M$.
\end{lemma}

\noindent{\em Proof:}  Fix a knot $D$ in $M$. Suppose $C$ represents an element of order $m$ in ${\rm H}_1(M;{\displaystyle\Bbb Z})$.
Since $lk_M(C,C')=0$, we may choose a two-sided surface $\Sigma$ with boundary $m$ times $C$
which is disjoint from a neighborhood of $C'$ and which meets $D$ transversely.  Then Dehn surgery
on $C'$ does not affect a neighborhood of $\Sigma$ and hence does not affect the oriented intersection of
$\Sigma$ and $D$.  The assertion follows.
$\Box$

We now show

\begin{lemma}
Let $C$ and $D$  be knots in a rational homology sphere $M$. Suppose $C$ is homologous to 0
and bounds a Seifert surface
$\Sigma$ disjoint from $D$ satisfying $lk_M(D,\alpha_j)=0$ for  $j=1,2,\ldots,2g$, where
$\alpha_1,\alpha_2,\ldots,\alpha_{2g}$ is a collection of simple closed curves on $\Sigma$ representing a basis
of ${\rm H}_1(\Sigma;{\displaystyle\Bbb Z})$.  
Let $M'$ denote the manifold resulting from $p/q$-Dehn surgery on $D$ 
in $M$, and let $C'$ denote the image of $C$ in $M'$.  Then  the Alexander 
polynomial of $C'$ in $M'$ is equal to the Alexander polynomial of $C$ in $M$.
\end{lemma}

\noindent
{\em Proof:} Let $\Sigma'$ be the image of $\Sigma$ in $M'$. 
Let $\alpha'_1, \alpha'_2, \ldots, \alpha'_{2g}$ denote the images of $\alpha_1, \alpha_2, \ldots, \alpha_{2g}$,
respectively, in $M'$.  The Alexander polynomial $\Delta_C$ of $C$ in $M$ is the 
determinant of the matrix $A$ with entries 
$a_{i,j} = lk_M(\alpha^+_i,\alpha_j) - t lk_M(\alpha^-_i,\alpha_j)$,
while the Alexander polynomial $\Delta_{C'}$ of $C'$ in $M'$ is the determinant 
of the matrix $A'$ with entries 
\mbox{$a'_{i,j} = lk_{M'}(\alpha^{'+}_i,\alpha'_j) - t lk_{M'}(\alpha^{'-}_i,\alpha'_j)$}.  
(Here $x^+$ and $x^-$ denote the plus- and minus-pushoffs of a simple closed curve $x$ on $\Sigma$
or $\Sigma'$.)  
(See [W,Appendix B], for example.)

But since $lk_M(D,\alpha_j)=0$ for $j=1,2,\ldots,2g$, we may apply Lemma A.5 to show that 
$ lk_{M'}(\alpha^{'+}_i,\alpha'_j) =  lk_M(\alpha^+_i,\alpha_j)$ and 
$lk_{M'}(\alpha^{'-}_i,\alpha'_j) =  lk_M(\alpha^-_i,\alpha_j)$ for all $i$ and $j$.  
It follows that $A' = A$, and hence $\Delta_{C'} = \Delta_C$.
$\Box$

\begin{prop}
Let $(K,L)$ be a boundary link in $S^3$, and let $\overline{L}$ be a lift 
of $L$ in the branched $k$-fold cover $\overline{S}^3_K$ for some $k$.  The 
Alexander polynomial of $\overline L$ in $\overline S^3_K$ is equal to the Alexander 
polynomial of $L$ in $S^3$.
\end{prop}

\noindent {\em Proof:}  Choose a link $C = (C_1,C_2,\ldots,C_n)$ in $S^3$ and 
surgery coefficients $r = (r_1,r_2,\ldots,r_n)$ so that $r$-Dehn surgery on $C$ 
yields $S^3$ and so that the image $K'$ of $K$ is an unknot. (See [R, Section  6D], 
for example.)  Isotoping $C$ as necessary, changing crossings of $C$ and $L$ but 
changing no crossings of $C \cup K$, we may assume that $L$, $K$, and $C_i$ bound Seifert surfaces
$\Sigma_L$, $\Sigma_K$, and $\Sigma_i$, respectively, with 
$\Sigma_L \cap \Sigma_K = \Sigma_L \cap \Sigma_i = \emptyset$ for 
$i = 1,2,\ldots,n$.

Let $L'$ denote the image of $L$ under $r$-surgery on $C$.  Let
$\alpha_1,\alpha_2,\ldots,\alpha_{2g}$ be a collection of simple closed curves on $\Sigma$ 
representing a basis
of ${\rm H}_1(\Sigma;{\displaystyle\Bbb Z})$.  By Lemma A.5, the linking number of the image of any component
$C_i$ with the image of any curve $\alpha_j$ is 0 at any stage of the sequence of surgeries on $C_1,C_2,\ldots,C_n$.  
Then we may apply
Lemma A.6 repeatedly to show that the Alexander 
polynomial $\Delta_{L'} = \Delta_L$.  Furthermore, 
the linking number of $K'$ and the image of any $\alpha_j$ in $S^3$ at the end of the 
surgery sequence is 0.  Then applying Lemma 3.2, we find that
$(L',K')$ and $(L',C'_i)$ are 
boundary links, where the $C'_i$ are the images of the $C_i$ in $S^3$ after 
the surgery.  In fact, since $K'$ is an unknot, we see that $(L',K')$ is splittable.

Now consider the $k$-fold branched cover of $S^3$ branched 
along the unknot $K'$, which is again $S^3$.  
Since $(L',K')$ is splittable, we see that  $L'$ has $k$ lifts 
$\overline L'_1,\overline L'_2,\ldots,\overline L'_k$ lying in disjoint 3-balls, 
and the Alexander polynomials satisfy $\Delta_{\overline L'_i} = \Delta_{L'}$ for 
$i = 1,2,\ldots,k$.  Let $\overline C'$ denote the inverse image of $C'$ in $\overline S^3_{K'}$.  

Finally, it is clear that surgery on $\overline C'$ with appropriate coefficients 
yields $\overline S^3_K$ and takes $\overline L'_i$ onto $\overline L_i$ for 
$i =1,2,\ldots,k$.  The assertion follows. $\Box$

\bibliographystyle{plain}

\begin{flushleft}

Department of Mathematics and Statistics
\newline The College of New Jersey
\newline Ewing, NJ 08628
\newline USA
\newline ccurtis@tcnj.edu
\end{flushleft}

\end{document}